\documentclass{article}%
\usepackage{amsmath}
\usepackage{amsfonts}
\usepackage{amssymb}
\usepackage{graphicx}%
\providecommand{\U}[1]{\protect\rule{.1in}{.1in}}

\newenvironment{proof}[1][Proof]{\noindent\textbf{#1.} }{\ \rule{0.5em}{0.5em}}
\begin{document}

\title{Eigenvalues and eigenfunctions of a boundary value problem }
\author{{\Large Erdo\u{g}an \c{S}en*}$^{1}${\Large , Azad Bayramov}*$^{2}$}
\date{}
\maketitle

{\Large *}$^{1}${\scriptsize {Department of Mathematics, Faculty of Arts and
Science, Nam\i k Kemal University, 59030, Tekirda\u{g}, Turkey..}}

{\scriptsize e-mail: esen@nku.edu.tr}$^{1},${\scriptsize abayramov@nku.edu.tr}%
$^{2}$



ABSTRACT. In this work a discontinuous boundary-value problem with retarded
argument which contains spectral parameter in the transmission conditions at
the point of discontinuity are investigated. We obtained asymptotic formulas
for the eigenvalues and eigenfunctions.

\section{Introduction}

Delay differential equations arise in many areas of mathematical modelling:
for example, population dynamics (taking into account the gestation times),
infectious diseases (accounting for the incubation periods), physiological and
pharmaceutical kinetics (modelling, for example, the body's reaction to
CO$_{2}$, etc. in circulating blood) and chemical kinetics (such as mixing
reactants), the navigational control of ships and aircraft, and more general
control problems.

Boundary value problems for differential equations of the second order with
retarded argument were studied in $[1-7]$, and various physical applications
of such problems can be found in $\left[  2\right]  $.

In the papers $[6,7]$ the asymptotic formulas for the eigenvalues and
eigenfunctions of discontinuous boundary value problem with retarded argument
and a spectral parameter in the boundary conditions were obtained.

The asymptotic formulas for the eigenvalues and eigenfunctions of
Sturm-Liouville problem with the spectral parameter in the boundary condition
were obtained in $\left[  8\right]  $.

The article $[9]$ is devoted to the study of the asymptotics of the solutions
to the Sturm-Liouville problem with the potential and the spectral parameter
having discontinuity of the first kind in the domain of definition of the solution.

In this paper we study the eigenvalues and eigenfunctions of discontinuous
boundary value problem with retarded argument and spectral parameters in the
transmission conditions. Namely we consider the boundary value problem for the
differential equation
\begin{equation}
y^{\prime\prime}(x)+q(x)y(x-\Delta(x))+\lambda y(x)=0 \tag*{(1)}%
\end{equation}
on $\left[  0,\frac{\pi}{2}\right)  \cup\left(  \frac{\pi}{2},\pi\right]  ,$
with boundary conditions%

\begin{equation}
y(0)\cos\alpha+y^{\prime}(0)\sin\alpha=0, \tag*{(2)}%
\end{equation}%
\begin{equation}
y(\pi)\cos\beta+y^{\prime}(\pi)\sin\beta=0, \tag*{(3)}%
\end{equation}
and transmission conditions%
\begin{equation}
y(\frac{\pi}{2}-0)-\sqrt[3]{\lambda}\delta y(\frac{\pi}{2}+0)=0, \tag*{(4)}%
\end{equation}%
\begin{equation}
y^{\prime}(\frac{\pi}{2}-0)-\sqrt[3]{\lambda}\delta y^{\prime}(\frac{\pi}%
{2}+0)=0, \tag*{(5)}%
\end{equation}
where the real-valued function $q(x)$ is continuous in $\left[  0,\frac{\pi
}{2}\right)  \cup\left(  \frac{\pi}{2},\pi\right]  ~$and has a finite limit
$q(\frac{\pi}{2}\pm0)=\lim_{x\rightarrow\frac{\pi}{2}\pm0}q(x),$ the real
valued function $\Delta(x)\geq0$ continuous in $\left[  0,\frac{\pi}%
{2}\right)  \cup\left(  \frac{\pi}{2},\pi\right]  $ and has a finite limit
$\Delta(\frac{\pi}{2}\pm0)=\lim_{x\rightarrow\frac{\pi}{2}\pm0}\Delta(x)$,
$x-\Delta(x)\geq0,$ if $\ x\in\left[  0,\frac{\pi}{2}\right)  ;x-\Delta
(x)\geq\frac{\pi}{2},$ if $x\in\left(  \frac{\pi}{2},\pi\right]  ;$ $\lambda$
is a real spectral parameter; $\delta$ is a arbitrary real number.

It must be noted that some problems with transmission conditions which arise
in mechanics (thermal condition problem for a thin laminated plate) were
studied in $\left[  10\right]  $.

Let $w_{1}(x,\lambda)$ be a solution of Eq. (1) on $\left[  0,\frac{\pi}%
{2}\right]  ,$ satisfying the initial conditions%
\begin{equation}
w_{1}\left(  0,\lambda\right)  =\sin\alpha,w_{1}^{\prime}\left(
0,\lambda\right)  =-\cos\alpha\tag*{(6)}%
\end{equation}
The conditions $(6)$ define a unique solution of Eq. $(1)$ on $\left[
0,\frac{\pi}{2}\right]  $($\left[  2\right]  $, p. 12).

After defining above solution we shall define the solution $w_{2}\left(
x,\lambda\right)  $ of Eq. (1) on $\left[  \frac{\pi}{2},\pi\right]  $ by
means of the solution $w_{1}\left(  x,\lambda\right)  $ by the initial
conditions%
\begin{equation}
w_{2}\left(  \frac{\pi}{2},\lambda\right)  =\lambda^{-1/3}\delta^{-1}%
w_{1}\left(  \frac{\pi}{2},\lambda\right)  ,\quad\omega_{2}^{\prime}%
({\frac{\pi}{2}},\>\lambda)=\lambda^{-1/3}\delta^{-1}\omega_{1}^{\prime
}({\frac{\pi}{2}},\>\lambda). \tag*{(7)}%
\end{equation}
The conditions (7) are defined as a unique solution of Eq. (1) on $\left[
\frac{\pi}{2},\pi\right]  .$

Consequently, the function $w\left(  x,\lambda\right)  $ is defined on
$\left[  0,\frac{\pi}{2}\right)  \cup\left(  \frac{\pi}{2},\pi\right]  $ by
the equality%
\[
w(x,\lambda)=\left\{
\begin{array}
[c]{ll}%
\omega_{1}(x,\lambda), & x\in\lbrack0,{\frac{\pi}{2}})\\
\omega_{2}(x,\lambda), & x\in({\frac{\pi}{2}},\pi]
\end{array}
\right.
\]
is a such solution of the Eq. (1) on $\left[  0,\frac{\pi}{2}\right)
\cup\left(  \frac{\pi}{2},\pi\right]  ;$which satisfies one of the. boundary
conditions and both transmission conditions.

\textbf{Lemma 1} \ Let $w\left(  x,\lambda\right)  $ be a solution of Eq.$(1)$
and\ $\lambda>0.$ Then the following integral equations hold:\
\begin{align}
w_{1}(x,\lambda)  &  =\sin\alpha\cos sx-\frac{\cos\alpha}{s}\sin
sx\tag*{(8)}\\
&  -\frac{1}{s}\int\limits_{0}^{{x}}q\left(  \tau\right)  \sin s\left(
x-\tau\right)  w_{1}\left(  \tau-\Delta\left(  \tau\right)  ,\lambda\right)
d\tau\text{ \ }\left(  s=\sqrt{\lambda},\lambda>0\right)  ,\nonumber
\end{align}%
\begin{align}
w_{2}(x,\lambda)  &  =\frac{1}{s^{2/3}\delta}w_{1}\left(  \frac{\pi}%
{2},\lambda\right)  \cos s\left(  x-\frac{\pi}{2}\right)  +\frac{w_{1}%
^{\prime}\left(  \frac{\pi}{2},\lambda\right)  }{s^{5/3}\delta}\sin s\left(
x-\frac{\pi}{2}\right) \tag*{(9)}\\
&  -\frac{1}{s}\int\limits_{\pi/2}^{{x}}q\left(  \tau\right)  \sin s\left(
x-\tau\right)  w_{2}\left(  \tau-\Delta\left(  \tau\right)  ,\lambda\right)
d\tau\text{ \ }\left(  s=\sqrt{\lambda},\lambda>0\right)  .\nonumber
\end{align}

\begin{proof}
To prove this, it is enough to substitute $\>-s^{2}\omega_{1}(\tau
,\lambda)-\omega_{1}^{\prime\prime}(\tau,\lambda)\>$ and $\>-s^{2}\omega
_{2}(\tau,\lambda)-\omega_{2}^{\prime\prime}(\tau,\lambda)\>$ instead of
$\>-q(\tau)\omega_{1}(\tau-\Delta(\tau),\lambda)\>$ and $\>-q(\tau)\omega
_{2}(\tau-\Delta(\tau),\lambda)\>$ in the integrals in (8) and (9)
respectively and integrate by parts twice.
\end{proof}

\textbf{Theorem 1 \ }The problem $(1)-(5)$ can have only simple eigenvalues.

\begin{proof}
Let $\widetilde{\lambda}$ be an eigenvalue of the problem $(1)-(5)$ and%
\[
\widetilde{y}(x,\widetilde{\lambda})=\left\{
\begin{array}
[c]{ll}%
\widetilde{y}_{1}(x,\widetilde{\lambda}), & x\in\lbrack0,{\frac{\pi}{2}}),\\
\widetilde{y}_{2}(x,\widetilde{\lambda}), & x\in({\frac{\pi}{2}},\pi]
\end{array}
\right.
\]
be a corresponding eigenfunction. Then from $(2)$ and $(6)$ it follows that
the determinant%
\[
W\left[  \widetilde{y}_{1}(0,\widetilde{\lambda}),w_{1}(0,\widetilde{\lambda
})\right]  =\left\vert
\begin{array}
[c]{c}%
\widetilde{y}_{1}(0,\widetilde{\lambda})\text{ \ \ \ \ \ \ \ \ \ }\sin\alpha\\
\widetilde{y}_{1}^{\prime}(0,\widetilde{\lambda})\text{ \ \ \ \ \ \ }%
-\cos\alpha
\end{array}
\right\vert =0,
\]
and by Theorem 2.2.2 in $\left[  2\right]  $ the functions $\widetilde{y}%
_{1}(x,\widetilde{\lambda})$ and $w_{1}(x,\widetilde{\lambda})$ are linearly
dependent on $\left[  0,\frac{\pi}{2}\right]  $. We can also prove that the
functions $\widetilde{y}_{2}(x,\widetilde{\lambda})$ and $w_{2}(x,\widetilde
{\lambda})$ are linearly dependent on $\left[  \frac{\pi}{2},\pi\right]  $.
Hence%
\begin{equation}
\widetilde{y}_{1}(x,\widetilde{\lambda})=K_{i}w_{i}(x,\widetilde{\lambda
})\text{ \ \ \ }\left(  i=1,2\right)  \tag*{(10)}%
\end{equation}
for some $K_{1}\neq0$ and $K_{2}\neq0$. We must show that $K_{1}=K_{2}$.
Suppose that $K_{1}\neq K_{2}$. From the equalities $(4)$ and $(10)$, we have%
\begin{align*}
\widetilde{y}(\frac{\pi}{2}-0,\widetilde{\lambda})-\sqrt[3]{\widetilde
{\lambda}}\delta\widetilde{y}(\frac{\pi}{2}+0,\widetilde{\lambda})  &
=\widetilde{y_{1}}(\frac{\pi}{2},\widetilde{\lambda})-\sqrt[3]{\widetilde
{\lambda}}\delta\widetilde{y_{2}}(\frac{\pi}{2},\widetilde{\lambda})\\
&  =K_{1}w_{1}(\frac{\pi}{2},\widetilde{\lambda})-\sqrt[3]{\widetilde{\lambda
}}\delta K_{2}w_{2}(\frac{\pi}{2},\widetilde{\lambda})\\
&  =\sqrt[3]{\widetilde{\lambda}}\delta K_{1}w_{2}(\frac{\pi}{2}%
,\widetilde{\lambda})-\sqrt[3]{\widetilde{\lambda}}\delta K_{2}w_{2}(\frac
{\pi}{2},\widetilde{\lambda})\\
&  =\sqrt[3]{\widetilde{\lambda}}\delta\left(  K_{1}-K_{2}\right)  w_{2}%
(\frac{\pi}{2},\widetilde{\lambda})=0.
\end{align*}
Since $\delta_{1}\left(  K_{1}-K_{2}\right)  \neq0$ it follows that
\begin{equation}
w_{2}\left(  \frac{\pi}{2},\widetilde{\lambda}\right)  =0. \tag{11}%
\end{equation}
By the same procedure from equality $(5)$ we can derive that%
\begin{equation}
w_{2}^{^{\prime}}\left(  \frac{\pi}{2},\widetilde{\lambda}\right)  =0.
\tag{12}%
\end{equation}
From the fact that $w_{2}(x,\widetilde{\lambda})$ is a solution of the
differential Eq. $(1)$ on $\left[  \frac{\pi}{2},\pi\right]  $ and satisfies
the initial conditions $(11)$ and $(12)$ it follows that $w_{1}(x,\widetilde
{\lambda})=0$ identically on $\left[  \frac{\pi}{2},\pi\right]  $ (cf. [2, p.
12, Theorem 1.2.1]).

By using we may also find%
\[
w_{1}\left(  \frac{\pi}{2},\widetilde{\lambda}\right)  =w_{1}^{^{\prime}%
}\left(  \frac{\pi}{2},\widetilde{\lambda}\right)  =0.
\]
From the latter discussions of $w_{2}(x,\widetilde{\lambda})$ it follows that
$w_{1}(x,\widetilde{\lambda})=0$ identically on $\left[  0,\frac{\pi}%
{2}\right)  \cup\left(  \frac{\pi}{2},\pi\right]  $. But this contradicts
$(6)$, thus completing the proof.
\end{proof}

\section{An existance theorem}

The function $\omega(x,\>\lambda)\>$defined in section $1$ is a nontrivial
solution of Eq. $(1)$ satisfying conditions $(2),(4)$ and $(5)$.
Putting$\>\omega(x,\>\lambda)\>$into $(3)$, we get the characteristic equation%
\begin{equation}
F(\lambda)\equiv\omega(\pi,\>\lambda)\cos\beta+\omega^{\prime}(\pi
,\>\lambda)\sin\beta=0. \tag*{(13)}%
\end{equation}

By Theorem 1.1 the set of eigenvalues of boundary-value problem (1)-(5)
coincides with the set of real roots of Eq. (13). Let $\>q_{1}=\int
\limits_{0}^{{\pi/2}}|q(\tau)|d\tau$ and $q_{2}=\int\limits_{{\pi/2}}^{{\pi}%
}\left\vert q(\tau)\right\vert d\tau.$

\textbf{Lemma 2} \ $(1)$ Let $\lambda\geq4q_{1}^{2}$. Then for the solution
$w_{1}\left(  x,\lambda\right)  $ of Eq. $(8)$, the following inequality
holds$:$%
\begin{equation}
\left\vert w_{1}\left(  x,\lambda\right)  \right\vert \leq\frac{1}{\left\vert
q_{1}\right\vert }\sqrt{4q_{1}^{2}\sin^{2}\alpha+\cos^{2}\alpha},\text{
\ \ }x\in\left[  0,\frac{\pi}{2}\right]  . \tag*{(14)}%
\end{equation}
$(2)$ Let $\lambda\geq\max\left\{  4q_{1}^{2},4q_{2}^{2}\right\}  $. Then for
the solution $w_{2}\left(  x,\lambda\right)  $ of Eq. $(9)$, the following
inequality holds$:$%
\begin{equation}
\left\vert w_{2}\left(  x,\lambda\right)  \right\vert \leq\frac{2\sqrt[3]{2}%
}{\sqrt[3]{q_{1}^{5}}\delta}\sqrt{4q_{1}^{2}\sin^{2}\alpha+\cos^{2}\alpha
},\text{\ \ }x\in\left[  \frac{\pi}{2},\pi\right]  \tag*{(15)}%
\end{equation}

\begin{proof}
Let $B_{1\lambda}=\max_{\left[  0,\frac{\pi}{2}\right]  }\left\vert
w_{1}\left(  x,\lambda\right)  \right\vert $. Then from $(8)$, it follows
that, for every $\lambda>0$, the following inequality holds:%
\[
B_{1\lambda}\leq\sqrt{\sin^{2}\alpha+\frac{\cos^{2}\alpha}{s^{2}}}+\frac{1}%
{s}B_{1\lambda}q_{1}.
\]
If $s\geq2q_{1}$ we get $(14)$. Differentiating $(8)$ with respect to $x$, we
have%
\begin{equation}
w_{1}^{\prime}(x,\lambda)=-s\sin\alpha\sin sx-\cos\alpha\cos sx-\int
\limits_{0}^{x}q(\tau)\cos s\left(  x-\tau\right)  w_{1}(\tau-\Delta\left(
\tau\right)  ,\lambda)d\tau\tag{16}%
\end{equation}
From $(16)$ and $(14)$, it follows that, for $s\geq2q_{1}$, the following
inequality holds:$.$%
\begin{equation}
\frac{\left\vert w_{1}^{\prime}(x,\lambda)\right\vert }{s^{5/3}}\leq\frac
{1}{\sqrt[3]{4q_{1}^{5}}}\sqrt{4q_{1}^{2}\sin^{2}\alpha+\cos^{2}\alpha}
\tag{17}%
\end{equation}
Let $B_{2\lambda}=\max_{\left[  \frac{\pi}{2},\pi\right]  }\left\vert
w_{2}\left(  x,\lambda\right)  \right\vert $. Then from $(9),(14)$ and $(17)$
it follows that, for $s\geq2q_{1}$ and $s\geq2q_{2}$, the following
inequalities holds:%
\begin{align*}
B_{2\lambda}  &  \leq\frac{2}{\sqrt[3]{4q_{1}^{5}}\delta}\sqrt{4q_{1}^{2}%
\sin^{2}\alpha+\cos^{2}\alpha}+\frac{1}{2q_{2}}B_{2\lambda}q_{2},\\
B_{2\lambda}  &  \leq\frac{2\sqrt[3]{2}}{\sqrt[3]{q_{1}^{5}}\delta}%
\sqrt{4q_{1}^{2}\sin^{2}\alpha+\cos^{2}\alpha}.
\end{align*}
Hence if $\lambda\geq\max\left\{  4q_{1}^{2},4q_{2}^{2}\right\}  $ we get
$(15)$.
\end{proof}

\textbf{Theorem 2} \ The problem $(1)-(5)$ has an infinite set of positive eigenvalues.

\begin{proof}
Differentiating $(9)$ with respect to$\>x$, we get%
\begin{align}
w_{2}^{\prime}(x,\lambda)  &  =-\frac{\sqrt[3]{s}}{\delta}w_{1}\left(
\frac{\pi}{2},\lambda\right)  \sin s\left(  x-\frac{\pi}{2}\right)
+\frac{w_{1}^{\prime}\left(  \frac{\pi}{2},\lambda\right)  }{\sqrt[3]{s^{2}%
}\delta}\cos s\left(  x-\frac{\pi}{2}\right) \nonumber\\
&  -\int\limits_{\pi/2}^{{x}}q(\tau)\cos s\left(  x-\tau\right)  w_{2}%
(\tau-\Delta\left(  \tau\right)  ,\lambda)d\tau.\text{ \ \ }(s=\sqrt{\lambda
},\lambda>0) \tag{18}%
\end{align}
From $(8),(9),(13),(16)$ and $(18)$, we get%
\begin{equation}
\left[  \frac{{1}}{s^{2/3}\delta}\biggl(\sin\alpha\cos{\frac{s\pi}{2}}%
-\frac{{\cos\alpha}}{s}\sin{\frac{s\pi}{2}}-{\frac{1}{s}}\int\limits_{0}%
^{{\frac{\pi}{2}}}q(\tau)\sin s({\frac{\pi}{2}}-\tau)\omega_{1}(\tau
-\Delta(\tau),\lambda)d\tau\biggr)\right. \nonumber
\end{equation}%
\[
\times\cos{\frac{s\pi}{2}}%
\]%
\[
-{\frac{1}{s^{5/3}\delta}}\biggl(s\sin\alpha\sin{\frac{s\pi}{2}}+\cos
\alpha\cos{\frac{s\pi}{2}}+\int\limits_{0}^{{\frac{\pi}{2}}}q(\tau)\cos
s({\frac{\pi}{2}}-\tau)\omega_{1}(\tau-\Delta(\tau),\lambda)d\tau\biggr)
\]%
\[
\left.  \times\sin{\frac{s\pi}{2}}-\frac{1}{s}\int\limits_{\pi/2}^{{\pi}%
}q(\tau)\sin s({\pi}-\tau)\omega_{2}(\tau-\Delta(\tau),\lambda)d\tau\right]
\cos\beta+
\]%
\[
\left[  -\frac{\sqrt[3]{s}}{\delta}\left(  \sin\alpha\cos{\frac{s\pi}{2}%
-}\frac{{\cos\alpha}}{s}\sin{\frac{s\pi}{2}-\frac{1}{s}}\int\limits_{0}%
^{{\frac{\pi}{2}}}q(\tau)\sin s({\frac{\pi}{2}}-\tau)\omega_{1}(\tau
-\Delta(\tau),\lambda)d\tau\right)  \right.
\]%
\[
\times\sin{\frac{s\pi}{2}}%
\]%
\[
-\frac{1}{\sqrt[3]{s^{2}}\delta}\left(  s\sin\alpha\sin{\frac{s\pi}{2}%
+\cos\alpha}\cos{\frac{s\pi}{2}+}\int\limits_{0}^{{\frac{\pi}{2}}}q(\tau)\cos
s({\frac{\pi}{2}}-\tau)\omega_{1}(\tau-\Delta(\tau),\lambda)d\tau\right)
\]%
\begin{equation}
\left.  \times\cos{\frac{s\pi}{2}}-\int\limits_{\frac{\pi}{2}}^{{\pi}}%
q(\tau)\cos s({\pi}-\tau)\omega_{2}(\tau-\Delta(\tau),\lambda)d\tau\right]
\sin\beta=0 \tag{19}%
\end{equation}
There are four possible cases:

$1.$ $\sin\alpha\neq0,\sin\beta\neq0;$

$2.$ $\sin\alpha\neq0,\sin\beta=0;$

$3.$ $\sin\alpha=0,\sin\beta\neq0;$

$4.$ $\sin\alpha=0,\sin\beta=0.$

In this paper we shall consider only case 1.The other cases may be considered
analogically. Let $\>\lambda\>$ be sufficiently large. Then, by $(14)$ and
$(15)$, the Eq. $(19)$ may be rewritten in the form%
\begin{equation}
\sqrt[3]{s}\sin s\pi+O(1)=0 \tag{20}%
\end{equation}
Obviously, for large$\>s\>$Eq. $(20)$ has an infinite set of roots. Thus the
theorem is proved.
\end{proof}

\section{Asymptotic Formulas for Eigenvalues and Eigenfunctions}

Now we begin to study asymptotic properties of eigenvalues and eigenfunctions.
In the following we shall assume that$\>s\>$is sufficiently large. From $(8)$
and $(14)$, we get
\[
\omega_{1}(x,\>\lambda)=O(1)\quad\mbox{on}\quad\lbrack0,\>{\frac{\pi}{2}%
}]\eqno(21)
\]
From $(9)$ and $(15)$, we get
\[
\omega_{2}(x,\>\lambda)=O(1)\quad\mbox{on}\quad\lbrack{\frac{\pi}{2}}%
,\>\pi]\eqno(22)
\]
The existence and continuity of the derivatives $\>\omega_{1s}^{\prime
}(x,\>\lambda)\>$for $\>0\leq x\leq{\frac{\pi}{2}},\>|\lambda|<\infty$,
and$\>\omega_{2s}^{\prime}(x,\>\lambda)\>$for $\>{\frac{\pi}{2}}\leq x\leq
\pi,\>|\lambda|<\infty$, follows from Theorem 1.4.1 in \cite{Nk}.

\textbf{Lemma 3} In case 1%

\begin{align}
\omega_{1s}^{\prime}(x,\>\lambda)  &  =O(1),\quad x\in\lbrack0,\>{\frac{\pi
}{2}}],\tag{23}\\
\omega_{2s}^{\prime}(x,\>\lambda)  &  =O(1),\quad x\in\lbrack\>{\frac{\pi}{2}%
},\>\pi]. \tag{24}%
\end{align}
are hold.

\begin{proof}
By differentiation of (1.8) with respect to $s$, we get, by (3.1)%
\begin{equation}
w_{1s}^{^{\prime}}(x,\lambda)=-\frac{1}{s}%
{\displaystyle\int\limits_{0}^{x}}
q(\tau)\sin s(x-\tau)w_{1s}^{^{\prime}}\left(  \tau-\Delta\left(  \tau\right)
,\lambda\right)  +Z(x,\lambda),\text{ \ }(\left\vert Z(x,\lambda)\right\vert
\leq Z_{0}) \tag{25}%
\end{equation}
Let $D_{\lambda}=\max_{[0,\frac{\pi}{2}]}\left\vert w_{1s}^{^{\prime}%
}(x,\lambda)\right\vert .$Then the existance of $D_{\lambda}$ follows from
continuity of derivation for $x\in\left[  0,\frac{\pi}{2}\right]  .$From (25)%
\[
D_{\lambda}\leq\frac{1}{s}q_{1}D_{\lambda}+Z_{0}.
\]
Now let $s\geq2q_{1}.$Then $D_{\lambda}\leq2Z_{0}$ and the validity of the
asymptotic formula (23) follows. Formula (24) may be proved analogically.
\end{proof}

\textbf{Theorem 3} \ Let $n$ be a natural number. For each sufficiently large
$n,$ in case $1,$ there is exactly one eigenvalue of the problem
$(1)-(5)$\ near $n^{2}.$

\begin{proof}
We consider the expression which is denoted by$\>O(1)$ in the Eq. $(20)$.%
\begin{align*}
&  \frac{\delta}{\sin\alpha\sin\beta}\left\{  -\frac{\sin(\alpha-\beta
)}{s^{2/3}\delta}\cos s\pi+\frac{\cos\alpha\cos\beta}{s^{5/3}\delta}\sin s\pi+%
{\displaystyle\int\limits_{0}^{\frac{\pi}{2}}}
\left[  \frac{\cos\beta}{s^{5/3}\delta}\sin s(\pi-\tau)+\right.  \right. \\
&  \left.  \frac{\sin\beta}{s^{2/3}\delta}\cos s\left(  \pi-\tau\right)
\right]  q\left(  \tau\right)  w_{1}\left(  \tau-\Delta\left(  \tau\right)
,\lambda\right)  d\tau+%
{\displaystyle\int\limits_{\frac{\pi}{2}}^{\pi}}
\left[  \frac{\cos\beta}{s}\sin s(\pi-\tau)+\sin\beta\cos s\left(  \pi
-\tau\right)  \right] \\
&  \left.  q\left(  \tau\right)  w_{2}\left(  \tau-\Delta\left(  \tau\right)
,\lambda\right)  d\tau\right\}
\end{align*}
If formulas $(21)-(23)$ are taken into consideration, it can be shown by
differentiation with respect to$\>s$ that for large$\>s\>$this expression has
bounded derivative.It is obvious that for large$\>s$ the roots of Eq. $(20)$
are situated close to entire numbers. We shall show that, for large$\>n$, only
one root $(20)$ lies near to each $n$. We consider the function $\>\phi
(s)=\sqrt[3]{s}\sin s\pi+O(1)\>$. Its derivative, which has the form$\>\phi
^{\prime}(s)=\frac{1}{3\sqrt[3]{s^{2}}}\sin s\pi+\sqrt[3]{s}\pi\cos\pi+O(1)$,
does not vanish for$\>s\>$close to$\>n\>$for sufficiently large$\>n\>$. Thus
our assertion follows by Rolle's Theorem.
\end{proof}

Let $\>n\>$ be sufficiently large. In what follows we shall denote
by$\>\lambda_{n}=s_{n}^{2}$ the eigenvalue of the problem $(1)-(5)$ situated
near $n^{2}$. We set $s_{n}=n+\delta_{n}$. From $(20)$ it follows that
$\delta_{n}=O\left(  \frac{1}{n^{1/3}}\right)  $. Consequently%
\[
s_{n}=n+O\bigl ({\frac{1}{n^{1/3}}}\bigr ).\eqno(26)
\]
The formula $(26)$ make it possible to obtain asymptotic expressions for
eigenfunction of the problem $(1)-(5)$.
From $(8),(16)$ and $(21)$, we get
\[
\omega_{1}(x,\>\lambda)=\sin\alpha\cos sx+O\bigl ({\frac{1}{s}}%
\bigr ),\eqno(27)
\]%
\[
\omega_{1}^{^{\prime}}(x,\>\lambda)=-s\sin\alpha\sin sx+O\bigl ({1}%
\bigr ).\eqno(28)
\]
From $(9),(22),(27)$ and $(28)$, we get%
\[
\omega_{2}(x,\>\lambda)={\frac{\sin\alpha}{s^{2/3}\delta}\cos}\frac{s\pi}%
{2}\cos s\left(  x-\frac{\pi}{2}\right)  -{\frac{\sin\alpha}{s^{2/3}\delta
}\sin}\frac{s\pi}{2}\sin s\left(  x-\frac{\pi}{2}\right)  +O\left(  \frac
{1}{s}\right)
\]%
\[
\omega_{2}(x,\>\lambda)={\frac{\sin\alpha}{s^{2/3}\delta}}\cos
sx+O\bigl ({\frac{1}{s}}\bigr )\eqno(29)
\]
By putting $(26)$ in the $(27)$ and $(29)$, we derive that%
\begin{align*}
u_{1n}  &  =w_{1}\left(  x,\lambda_{n}\right)  =\sin\alpha\cos
nx+O\bigl ({\frac{1}{n^{1/3}}}\bigr ),\\
u_{2n}  &  =w_{2}\left(  x,\lambda_{n}\right)  ={\frac{{\sin\alpha}}{\delta
n^{2/3}}}\cos nx+O\bigl ({\frac{1}{n}}\bigr ).
\end{align*}
Hence the eigenfunctions$\>u_{n}(x)\>$have the following asymptotic
representation:%
\[
u_{n}(x)=\left\{
\begin{array}
[c]{lll}%
\sin\alpha\cos nx+O\bigl ({\frac{1}{n^{1/3}}}\bigr ), & \mbox{for} &
x\in\lbrack0,{\frac{\pi}{2}}),\\
{\frac{{\sin\alpha}}{\delta n^{2/3}}}\cos nx+O\bigl ({\frac{1}{n}}\bigr ) &
\mbox{for} & x\in({\frac{\pi}{2}},\pi].
\end{array}
\right.
\]
Under some additional conditions the more exact asymptotic formulas which
depend upon the retardation may be obtained. Let us assume that the following
conditions are fulfilled:

\noindent\textbf{a)} The derivatives $\>q^{\prime}(x)\>$ and $\>\Delta
^{\prime\prime}(x)\>$ exist and are bounded in$\>[0,{\frac{\pi}{2}}%
)\bigcup({\frac{\pi}{2}},\pi]\>$ and have finite limits $\>q^{\prime}%
({\frac{\pi}{2}}\pm0)=\lim\limits_{x\rightarrow{\frac{\pi}{2}}\pm0}q^{\prime
}(x)\>$ and $\>\Delta^{\prime\prime}({\frac{\pi}{2}}\pm0)=\lim
\limits_{x\rightarrow{\frac{\pi}{2}}\pm0}\Delta^{\prime\prime}(x)$, respectively.

\noindent\textbf{b)} $\>\Delta^{\prime}(x)\leq1 \>$ in$\>[0,{\frac{\pi}{2}%
})\bigcup({\frac{\pi}{2}},\pi]$, $\>\Delta(0)=0\>$ and $\>\lim
\limits_{x\rightarrow{\frac{\pi}{2}}+ 0}\Delta(x)=0$.

By using b), we have
\begin{align}
x-\Delta(x)  &  \geq0,\>x\in\lbrack0,{\frac{\pi}{2}})\tag{30}\\
x-\Delta(x)  &  \geq{\frac{\pi}{2}},\>x\in({\frac{\pi}{2}},\pi] \tag{31}%
\end{align}

From $(27)$, $(29)$, $(30)$ and $(31)$ we have%
\begin{equation}
w_{1}\left(  \tau-\Delta\left(  \tau\right)  ,\lambda\right)  =\sin\alpha\cos
s\left(  \tau-\Delta\left(  \tau\right)  \right)  +O\bigl ({\frac{1}{s}%
}\bigr ), \tag{32}%
\end{equation}%
\begin{equation}
w_{2}\left(  \tau-\Delta\left(  \tau\right)  ,\lambda\right)  ={\frac
{{\sin\alpha}}{s^{2/3}\delta}}\cos s\left(  \tau-\Delta\left(  \tau\right)
\right)  +O\bigl ({\frac{1}{s}}\bigr ). \tag{33}%
\end{equation}
Putting these expressions into $(19)$, we have%
\begin{align}
0  &  =-\frac{s^{1/3}}{\delta}\sin\alpha\sin\beta\sin s\pi+\frac{\sin\left(
\alpha-\beta\right)  }{s^{2/3}\delta}\cos s\pi-\frac{\sin\alpha\sin\beta
}{s^{2/3}\delta}\nonumber\\
&  \times\left\{  \cos s\pi%
{\displaystyle\int\limits_{0}^{\pi}}
\frac{q\left(  \tau\right)  }{2}\left[  \cos s\Delta(\tau)+\cos s\left(
2\tau-\Delta(\tau)\right)  \right]  d\tau\right. \nonumber\\
&  \left.  +\sin s\pi%
{\displaystyle\int\limits_{0}^{\pi}}
\frac{q\left(  \tau\right)  }{2}\left[  \sin s\Delta(\tau)+\sin s\left(
2\tau-\Delta(\tau)\right)  \right]  d\tau\right\}  +O\left(  \frac{1}{s^{5/3}%
}\right)  \tag{34}%
\end{align}

\noindent Let%
\begin{align}
K\left(  x,s,\Delta(\tau)\right)   &  ={\frac{1}{2}}\int\limits_{0}^{x}%
q(\tau)\sin s\Delta(\tau)d\tau,\nonumber\\
L(x,s,\Delta(\tau))  &  ={\frac{1}{2}}\int\limits_{0}^{x}q(\tau)\cos
s\Delta(\tau)d\tau\tag{35}%
\end{align}
\noindent It is obviously that these functions are bounded for $\>0\leq
x\leq\pi,\>0<s<+\infty$.

\noindent Under the conditions a) and b) the following formulas%
\begin{equation}
\int\limits_{0}^{x}q(\tau)\cos s(2\tau-\Delta(\tau))d\tau=O\left(  {\frac
{1}{s}}\right)  ,\quad\int\limits_{0}^{x}q(\tau)\sin s(2\tau-\Delta
(\tau))d\tau=O\left(  {\frac{1}{s}}\right)  \tag{36}%
\end{equation}
\noindent can be proved by the same technique in Lemma 3.3.3 in \cite{Nk}.
From $(34),(35)$ and$\>(36)\>$we have%
\[
\sin s\pi\left[  s\sin\alpha\sin\beta+K\left(  \pi,s,\Delta\left(
\tau\right)  \right)  \sin\alpha\sin\beta\right]  -
\]%
\[
\cos s\pi\left[  \sin\alpha\cos\beta-\cos\alpha\sin\beta-L\left(  \pi
,s,\Delta\left(  \tau\right)  \right)  \sin\alpha\sin\beta\right]  +O\left(
\frac{1}{s}\right)  =0.
\]
Hence%
\[
\tan s\pi=\frac{1}{s}\left[  \cot\beta-\cot\alpha-L\left(  \pi,s,\Delta\left(
\tau\right)  \right)  \right]  +O\left(  \frac{1}{s^{2}}\right)  .
\]
Again if we take $s_{n}=n+\delta_{n}$, then
\[
\tan\left(  n+\delta_{n}\right)  \pi=\tan\delta_{n}\pi=\frac{1}{n}\left[
\cot\beta-\cot\alpha-L\left(  \pi,n,\Delta\left(  \tau\right)  \right)
\right]  +O\left(  \frac{1}{n^{2}}\right)
\]
hence for large $n$,%
\[
\delta_{n}=\frac{1}{n\pi}\left[  \cot\beta-\cot\alpha-L\left(  \pi
,n,\Delta\left(  \tau\right)  \right)  \right]  +O\left(  \frac{1}{n^{2}%
}\right)
\]

\noindent and finally%
\begin{equation}
s_{n}=n+\frac{1}{n\pi}\left[  \cot\beta-\cot\alpha-L\left(  \pi,n,\Delta
\left(  \tau\right)  \right)  \right]  +O\left(  \frac{1}{n^{2}}\right)
\tag{37}%
\end{equation}
Thus, we have proven the following theorem.

\textbf{Theorem 4} \ If conditions a) and b) are satisfied then, the positive
eigenvalues $\lambda_{n}=s_{n}^{2}\>$ of the problem (1)-(5) have the $(37)$
asymptotic representation for$\>n\rightarrow\infty\>$.

We now may obtain a sharper asymptotic formula for the eigenfunctions. From
$(8)$ and $(32)$%
\begin{align*}
w_{1}(x,\lambda)  &  =\sin\alpha\cos sx-\frac{\cos\alpha}{s}\sin sx\\
&  -\frac{\sin\alpha}{s}\int\limits_{0}^{{x}}q\left(  \tau\right)  \sin
s\left(  x-\tau\right)  \cos s\left(  \tau-\Delta\left(  \tau\right)  \right)
d\tau+O\left(  \frac{1}{s^{2}}\right)  .
\end{align*}
Thus, from $(35)$ and $(36)$%
\[
w_{1}(x,\lambda)=\sin\alpha\cos sx\left[  1+\frac{K\left(  x,\>s,\>\Delta
(\tau)\right)  }{s}\right]
\]%
\begin{equation}
-\frac{\sin sx}{s}\left[  \cos\alpha+\sin\alpha L\left(  x,\>s,\>\Delta
(\tau)\right)  \right]  +O\left(  \frac{1}{s^{2}}\right)  . \tag{38}%
\end{equation}
\noindent Replacing$\>s$ by$\>s_{n}$ and using $(37)$, we have
\[
u_{1n}(x)=w_{1}(x,\lambda_{n})=\sin\alpha\left\{  \cos nx\left[
1+\frac{K\left(  x,n,\Delta\left(  \tau\right)  \right)  }{n}\right]  \right.
\]%
\begin{equation}
\left.  -\frac{\sin nx}{n\pi}\left[  \left(  \cot\beta-\cot\alpha-L\left(
\pi,n,\Delta\left(  \tau\right)  \right)  \right)  x+\left(  \cot
\alpha+L\left(  x,n,\Delta\left(  \tau\right)  \right)  \right)  \pi\right]
\right\}  +O\left(  \frac{1}{n^{2}}\right)  \tag{39}%
\end{equation}
From $(16),(32)$ and $(35)$, we have%
\begin{align}
\frac{w_{1}^{^{\prime}}\left(  x,\lambda\right)  }{s^{5/3}}  &  =-\frac
{\sin\alpha\sin sx}{s^{2/3}}\left(  1+\frac{K\left(  x,\>s,\>\Delta
(\tau)\right)  }{s}\right) \nonumber\\
&  -\frac{\cos sx}{s^{5/3}}\left(  \cos\alpha+\sin\alpha L\left(
x,\>s,\>\Delta(\tau)\right)  \right)  +O\left(  \frac{1}{s^{2}}\right)  \text{
},x\in\left(  0,\frac{\pi}{2}\right]  \tag{40}%
\end{align}
From $(9),(33),(36),(38)$ and $(40),$ we have%
\begin{align*}
w_{2}\left(  x,\lambda\right)   &  =\frac{1}{s^{2/3}\delta}\left\{  \sin
\alpha\cos\frac{s\pi}{2}\left[  1+\frac{K\left(  \frac{\pi}{2},s,\Delta
(\tau)\right)  }{s}\right]  -\frac{\sin\frac{s\pi}{2}}{s}\left[  \cos
\alpha+\sin\alpha L\left(  \frac{\pi}{2},s,\Delta(\tau)\right)  \right]
\right. \\
&  \left.  +O\left(  \frac{1}{s^{2}}\right)  \right\}  \cos s\left(
x-\frac{\pi}{2}\right)  -\frac{1}{\delta}\left\{  \sin\alpha\sin\frac{s\pi}%
{2}\left[  1+\frac{K\left(  \frac{\pi}{2},s,\Delta(\tau)\right)  }{s}\right]
\right. \\
&  \left.  -\frac{\cos\frac{s\pi}{2}}{s}\left[  \cos\alpha+\sin\alpha L\left(
\frac{\pi}{2},s,\Delta(\tau)\right)  \right]  +O\left(  \frac{1}{s^{2}%
}\right)  \right\}  \sin s\left(  x-\frac{\pi}{2}\right) \\
&  -\frac{1}{s}\int\limits_{\pi/2}^{{x}}q\left(  \tau\right)  \sin s\left(
x-\tau\right)  \left[  \frac{\sin\alpha}{s^{2/3}\delta}\cos s\left(
\tau-\Delta\left(  \tau\right)  \right)  +O\left(  \frac{1}{s}\right)
\right]  d\tau\\
&  =\frac{\sin\alpha}{s^{2/3}\delta}\cos sx\left[  1+\frac{K\left(  \frac{\pi
}{2},s,\Delta(\tau)\right)  }{s}\right]  -\frac{\sin sx}{s^{5/3}\delta}\left(
\cos\alpha+\sin\alpha L\left(  \frac{\pi}{2},s,\Delta(\tau)\right)  \right) \\
&  -\frac{\sin\alpha}{s^{5/3}\delta}\int\limits_{\pi/2}^{{x}}\frac{q\left(
\tau\right)  }{2}\left[  \sin s\left(  x-\Delta(\tau)\right)  +\sin s\left(
x-\left(  2\tau-\Delta(\tau)\right)  \right)  \right]  d\tau+O\left(  \frac
{1}{s^{2}}\right) \\
&  =\frac{\sin\alpha}{s^{2/3}\delta}\cos sx\left[  1+\frac{K\left(
x,s,\Delta(\tau)\right)  }{s}\right] \\
&  -\frac{\sin sx}{s^{5/3}\delta}\left(  \cos\alpha+\sin\alpha L\left(
x,s,\Delta(\tau)\right)  \right)  +O\left(  \frac{1}{s^{2}}\right)  ,\text{
\ \ }x\in\left(  \frac{\pi}{2},\pi\right]  .
\end{align*}
\noindent Now, replacing $\>s\>$ by $\>s_{n}\>$and using $(37)$, we have
\[
u_{2n}(x)=\frac{\sin\alpha}{n^{2/3}\delta}\left\{  \cos nx\left[
1+\frac{K\left(  x,n,\Delta\left(  \tau\right)  \right)  }{n}\right]
-\frac{\sin nx}{n^{5/3}\pi}\right.
\]%
\[
\left.  \times\left[  \left(  \cot\beta-\cot\alpha-L\left(  \pi,n,\Delta
\left(  \tau\right)  \right)  \right)  x+\left(  \cot\alpha+L\left(
x,n,\Delta\left(  \tau\right)  \right)  \right)  \pi\right]  \right\}
+O\bigl({\frac{1}{{n^{2}}}}\bigr
).\eqno(41)
\]
Thus, we have proven the following theorem.

\textbf{Theorem 5} \ If conditions a) and b) are satisfied then, the
eigenfunctions$\>u_{n}(x)\>$of the problem (1)-(5) have the following
asymptotic representation for$\>n\rightarrow\infty\>$:
\[
u_{n}(x)=\left\{
\begin{array}
[c]{lll}%
u_{1n}(x) & \mbox{for} & x\in\lbrack0,{\frac{\pi}{2}})\\
&  & \\
u_{2n}(x) & \mbox{for} & x\in({\frac{\pi}{2}},\pi]
\end{array}
\right.
\]
\noindent where$\>u_{1n}(x)\>$and$\>u_{2n}(x)\>$defined as in (39) and (41) respectively.

\section{Conclusion}

In this study firstly we obtain asymptotic formulas for eigenvalues and
eigenfunctions for discontinuous boundary value problem with retarded argument
which contains spectral parameter in the transmission conditions. Then under
additional conditions a.) and b.) the more exact asymptotic formulas, which
depend upon the retardation obtained.

\end{document}